\documentclass[12pt]{amsart}
\usepackage{mathtools, amssymb, bm, amsthm, amsmath, mathrsfs}

\usepackage{geometry}
\usepackage{biblatex}
\usepackage{color}
\usepackage{amsthm}
\usepackage{lmodern}

\usepackage{hyperref}

\newtheorem{theorem}{Theorem}

\newtheorem{mainthm}{Theorem}

\newtheorem{proposition}[theorem]{Proposition}
\newtheorem{lemma}[theorem]{Lemma}

\newtheorem{conjecture}[theorem]{Conjecture}

\newtheoremstyle{named}{}{}{\itshape}{}{\bfseries}{.}{.5em}{\thmnote{#3's }#1}
\theoremstyle{named}

\theoremstyle{definition}
\newtheorem{remark}[theorem]{Remark}

\numberwithin{theorem}{section}
\numberwithin{equation}{section}

\newcommand{\R}{\mathbb R}
\newcommand{\Z}{\mathbb Z}
\newcommand{\N}{\mathbb N}
\newcommand{\Q}{\mathbb Q}

\newcommand{\eps}{\varepsilon}

\newcommand{\sqfull}{S^{{\mathstrut\hspace{0.05em}\blacksquare}}_f}
\newcommand{\sqfullid}{S^{{\mathstrut\hspace{0.05em}\blacksquare}}_{\text{Id}}}

\newcommand{\sqfree}{S_f^{\mathstrut\hspace{0.05em}\square}}

\newcommand{\vecg}{\mathbf{g}}

\begin{document}

\title{Square-full values of quadratic polynomials}

\author{Watcharakiete Wongcharoenbhorn}

\author{Yotsanan Meemark}

\address{Department of Mathematics and Computer Science, Faculty of Science, Chulalongkorn University, Bangkok, Thailand 10330}

\email{w.wongcharoenbhorn@gmail.com}
\email{yotsanan.m@chula.ac.th}


\maketitle

\begin{abstract}
A \textit{square-full} number is a positive integer for which all its prime divisors divide itself at least twice. We have known that for a relatively prime pair $(a,b)\in\N\times \N\cup\{0\}$ with an affine polynomial $f(x)=ax+b$, the number of $n\leqslant N$ such that $f(n)$ is square-full is $\asymp_{a,b} N^\frac{1}{2}$. For $\eps>0$ and $f(x)\in\Z[x]$ is an admissible quadratic polynomial, we show that the number of square-full integers assuming $f(n)$ for $n\leqslant N$ is at most
$$O_{\eps,f}(N^{\varpi+\eps})$$ 
for some absolute constant $\varpi<1/2$. Under the assumption on the $abc$ conjecture, we expect the upper bound to be $O_{\eps,f}(N^\varepsilon)$.
\end{abstract}

\section{Introduction}
\label{sec:introduction}
A positive integer $n$ is called \textit{square-full} if for every prime $p$ dividing $n$, then $p^2\mid n$. We may put $1$ to be a square-full number by this definition. Throughout this work, $f(x)$ is a polynomial with integer coefficients. Let us denote $\sqfull(N)$ the number of positive integers $n\leqslant N$ such that $f(n)$ is square-full. It is known for an identity polynomial $\text{Id}(n)=n$ that
\begin{align*}
    \sqfullid(N)=\dfrac{\zeta(3/2)}{\zeta(3)}N^{\frac{1}{2}}+\dfrac{\zeta(2/3)}{\zeta(2)}N^{\frac{1}{3}}+o(N^{\frac{1}{6}}).
\end{align*}
This result is due to P. T. Bateman and E. Grosswald \cite{Bateman-Grosswald} for which they improved the pioneer result \cite{square-full-Erdos-Szekeres} of P. Erd\H{o}s and G. Szekeres by producing the term of order $N^{\frac{1}{3}}$ explicitly. They also indicate in the paper that the exponent $1/6$ in the error term cannot be improved without proving some kind of the quasi-Riemann hypothesis. Furthermore, there are studies on $f(n)=an+b$ for $(a,b)\not=(1,0)$ as well, and the counting function still is of order $N^{\frac{1}{2}}$ up to a constant depending only on $a,b$. For instance, T. H. Chan \cite{square-full-inAP} provided the asymptotic formula for such a polynomial with $\gcd(a,b)=1$, and improved the error term for large $a$ from previous works in this direction.\\

In 1931, T. Estermann \cite{square-free-Estermann} studied the counting function for quadratic polynomials assuming square-free values. Here, \textit{square-free} integer is similarly defined, for which we replace $p^2\mid n$ by $p\|n$ in our definition of square-full numbers. Analogous to the square-full notation, we denote $\sqfree(N)$ by the number of positive integers $n\leqslant N$ such that $f(n)$ is square-free. Particularly, he proved for $f(x)=x^2+h$ with a nonzero integer $h$ that 
$$\sqfree(N)=\prod_{p}\left(1-\dfrac{\rho_h(p^2)}{p^2}\right) N+O(N^{\frac{2}{3}}\log N),$$
where $\rho_h(m):=\#\{n\in\Z/ m\Z: n^2+h\equiv 0\mod{m}\}$ for all $m\geqslant 2$. After that, around 80 years later, D. R. Heath-Brown \cite{square-free} improved in 2012 the error term for $h=1$ of the above formula to be $O_\varepsilon(N^{\frac{7}{12}+\varepsilon})$ for arbitrarily small $\varepsilon>0$.\\

Our interest is to determine the number of square-full numbers of quadratic form because it seems possible to apply the method of Heath-Brown used in \cite{square-free}. It turns out that there might be bias toward it, as suggested by our first main result on special quadratic polynomials.
\begin{mainthm}
%\label\ref{thm:squarefull-on-quadratic-forms}
\label{thm:squarefull-on-quadratic-forms}{}
Let $\eps>0$ be given and $N$ be a large number depending on $\eps$. Let $\alpha\in\N$ and $f(x)=x^2+\alpha^2$. The number of $n\leqslant N$ for which $f(n)$ is square-full is, uniformly for $\alpha\leqslant N^\varpi$, at most $$O_{\eps}\left(N^{\varpi+\eps}\right),$$
where $\varpi=\frac{29}{100}$.
\end{mainthm}

Here, we employ the celebrated determinant method introduced by Heath-Brown. His key idea is to leverage properties of the sum of two squares and the unique factorization of $\Z[i]$ to extract an equation involving a bi-homogeneous polynomial. Then, he constructed another equation with another bi-homogeneous form, in which it vanishes at each solution of the previous equation, to help control the count of solutions. \\

We are also interested in the more general case of a quadratic polynomial $f(x)=ax^2+bx+c\in \Z[x]$. An important property of any square-full integer is that it can be written uniquely as $e^2d^3$ for an integer $e$ and a square-free $d$. It is seen that $\sqfull(N)\leqslant S^{{\mathstrut\hspace{0.05em}\blacksquare}}_g(2aN)$ if $g(x)=px^2+q$ with $p=a$ and $q=a(b^2-4ac)$ since
\begin{align}
\label{eq:majorant-g}
    ax^2+bx+c=\frac{p(2ax+b)^2+q}{4a^2}
\end{align}
and the left-hand side is of the form $e^2d^3$ if $py^2+q=(2ae)^2d^3=:e'^2d^3$. We have that $$S^{{\mathstrut\hspace{0.05em}\blacksquare}}_f(N)\leqslant \#\{n\leqslant 2aN: e'^2d^3=pn^2+q, d\text{ square-free}\}=S^{{\mathstrut\hspace{0.05em}\blacksquare}}_g(2aN).$$ We call this $g$ a \textit{majorant of $f$}. Hence, it suffices to consider $f(x)=px^2+q$ for general polynomials. We are particularly interested in polynomials for which their majorant $g(x)=px^2+q$ satisfies $q\not=0$, i.e., when its discriminant is nonzero. Let us call such a polynomial $f(x)$ \textit{admissible} throughout the paper.
\\

To formulate the general result, we define 
\begin{align}
\label{eq:varpi_0-def}
    \varpi_0=\frac{2\beta}{3\log 3}\approx 0.1688,
\end{align}
where $\beta:=\frac{2+\sqrt{3}}{2\sqrt{3}}\log \frac{2+\sqrt{3}}{2\sqrt{3}}-\frac{2-\sqrt{3}}{2\sqrt{3}}\log \frac{2-\sqrt{3}}{2\sqrt{3}}$. We now state our second main theorem, which establishes the bound in general that remains strictly less than the $\frac{1}{2}$-power of $N$. The main ingredient is the improved bound of Ellenberg-Venkatesh \cite{Ellenberg-Venkatesh} on counting $3$-torsion of the class group $\Q(\sqrt{-D})$ with the Helfgott-Venkatesh mechanism \cite{Helfgott-Venkatesh}. These provide the bound on counting integral points on Mordell curves.
\begin{mainthm}
\label{thm:general-case}
For $\eps>0$ and an admissible quadratic polynomial $f(x)$, the number of $n\leqslant N$ for which $f(n)$ is square-full is at most $$O_{\eps,f}\left(N^{\varpi+\eps}\right)$$
where $\varpi=\frac{2(1+4\varpi_0)}{5+12\varpi_0}\approx0.4769$. Moreover, if the Mordell curve conjecture (Conjecture \ref{conj:weak-effective-Mordell} below) is true, then we can replace $\varpi$ by $\frac{2}{5}$.
\end{mainthm}

Moreover, we note that the problem of counting square-full values of a certain quadratic polynomial such as $f(x)=x^2+h$ for a nonzero integer $h$, is related to counting pairs of quasi-consecutive square-full numbers, since $n^2$ itself is square-full. In fact, if we fix a polynomial $f(x)=x(x+1)$, then counting pairs of consecutive square-full integers that do not exceed $N$ is the same as determining $\sqfull(N)$. T. Reuss \cite{consec-square-full} studied this problem, as well as some related problems on $k$-free numbers. Our result for special $f(x)$ is the same as his result, which also shows that the number of such pairs of integers up to $N$ is $O_\eps(N^{\frac{29}{100}+\eps})$. There are 2 challenges in adapting the determinant method in our problem. First, we need to deal with more than one nonlinear term during the extraction from the sum of two squares to exploit the method efficiently. We are helped by the fact that a square-full integer can be written as $e^2d^3$ for a square-free integer $d$ uniquely, and the exponent $2$ allows us to revert to a single linear term via Taylor's theorem. Second, we need to modify the counting argument for a certain subcase where we do not have as effective a result on elliptic curves as Thue's equations. We modify to counting instead on Thue's equation, which is possible by our new vanishing equation, leveraging the effective result of Bombieri-Schmidt \cite{bombieri-schmidt}.\\ 

Recently, T. Browning and I. Shparlinski \cite{square-free-random} have proved that almost all polynomials of arbitrary degree $k\geqslant 2$, in the sense of naive height, take infinitely many square-free values. We give the final remark that almost all quadratic polynomials, in this sense, have small bounds for their square-full counting formulae (see Remark \ref{rm:rand-polynom}).

\section{Preliminaries}
\label{sec:Prelim}
We show that, although it seems difficult to find a positive integer $n$ so that $n^2+4$ is square-full, there are, in fact, infinitely many such values of $n$. We may let $n=2d$ be an even number, and rewrite $n^2+4=4(d^2+1)$. There are infinitely many solutions $(d,k)\in\Z^2$ such that
$$4(d^2+1)=8k^2\Longleftrightarrow d^2-2k^2=-1,$$
since it is a negative Pell's equation with $(1,1)$ as a base solution. Since $8k^2$ is always square-full, we produce infinitely many square-full numbers of the form $n^2+4$. Also, we shall note on a polynomial $f(x)$ of degree at least three with no repeated roots that it seems impossible that $f(n)$ will assume infinite values of square-full integers. This is thanks to the result of A. Granville \cite{Granville} as a consequence of the $abc$ conjecture; see Theorem 7 of \cite{Granville} and the discussion around it. Our example $f(x)=x^2+4$ shows that the number of $n\leqslant N$ for which $f(n)$ is square-full is $\gg\log N$. We may ask, as mentioned by Granville, on counting consecutive square-full numbers that
$$\sqfull(N)\sim c_f\log N,$$
for some constant $c_f$. In particular, this might hold for any admissible quadratic polynomial $f(x)$ with no fixed prime $p$ such that $p\|f(n)$ for all sufficiently large $n$. The constant $c_f$ is allowed to depend on the coefficients of $f(x)$. In Proposition \ref{prop:general-polyn-via-abc}, we shall see that assuming the $abc$ conjecture, the left-hand side is bounded by $C_{\eps,f} N^\eps$ for some constant $C_{\eps, f}>0$. This result is perhaps far from optimal. It is therefore interesting if its asymptotic behavior does exist. By the Erd\H{o}s-Szekeres theorem, we have that the probability of $n$ in each dyadic range $[\frac{N}{2},N]$ to be a square-full is of order $1/\sqrt{N/2}$. Then, we expect the probability of both $n,n+1$ in the same range to be square-full to be of order $1/(N/2)$. Thus, it is reasonable to expect the number of consecutive square-full integers up to $N$ to be of order
$$\dfrac{N}{2}\cdot\dfrac{1}{N/2}+\dfrac{N}{4}\cdot\dfrac{1}{N/4}+\cdots \ll\log N,$$
where the sum ranges through all dyadic intervals not exceeding $N$. This gives a simple heuristic argument that supports the aforementioned question.
\\

For a fixed polynomial $f(x)$, it is convenient to study the expression
\begin{equation}
\label{eq:split-into-dyadic-of-counting-sqfull}
\sqfull(2N)-\sqfull(N) =\sum_{\substack{N<n\leqslant 2N\\ e^2d^3=f(n)}} \mu^2(d).
\end{equation}
This is because finding the upper bound for it implies the bound for $\sqfull(N)$ itself by a dyadic sum, at the cost of a logarithmic factor. By defining $\mathcal{M}(E,D) := \#\{(e,d,n)\in\N^3: N<n\leqslant 2N, E<e\leqslant 2E, D<d\leqslant 2D, e^2d^3=f(n)\}$, the right-hand side of Equation \eqref{eq:split-into-dyadic-of-counting-sqfull} is at most
\begin{align*}
\sum_{\ell\ll \log N} \sum_{\substack{\frac{N}{2^{\ell}}<E\leqslant \frac{N}{2^{\ell-1}}\\E^2D^3\asymp_f N^2}}\mathcal{M}(E,D)\ll (\log N)\sup_{\substack{E\leqslant N\\ E^2D^3\asymp_f N^2}}\mathcal{M}(E,D),
\end{align*}
where the sums are dyadic sums over $E$. From this, the problem is reduced to obtaining the bound for $\mathcal{M}(E,D)$ in the whole range. In his paper \cite{square-free-Estermann}, Estermann showed that the contribution to $\mathcal{M}(E,D)$ for fixed $d\in [D,2D]\subset [1,N]$ is $O(\log N)$ so we have $\mathcal{M}(E,D)\ll D\log N$.\\

It seems difficult to extend Heath-Brown's determinant method to a polynomial $f(x)=x^2+h$ with a non-perfect square $h$. For general quadratic polynomials, we use a different approach by employing the bound on counting integral points on Mordell curves, which is defined as $y^2=x^3+D$ for a nonzero integer $D$. This general case will be carried out in the following Section \ref{sec:Mordell-curves-counting}. By convention, for a positive integer $n$, we denote $\text{rad}(n)=\prod_{p|n} p$. This notation will be used in stating the $abc$ conjecture (Conjecture \ref{conj:abc-conjecture} below) to prove the conditional result for this general case.\\

In Section \ref{sec:determinant-method-and-special-polyn}, we prove our first main theorem on a special polynomial $f(x)=x^2+\alpha^2$. We proceed in parallel to Heath-Brown's work \cite{square-free} on square-free values of $n^2+1$, by refining the bound of $\mathcal{M}(E,D)$ to be $N^\eta\sqrt{M}\min(\sqrt{E},\sqrt{D})$, with a new parameter $M$. The idea is to construct, using this $M$, another bi-homogeneous form for which all the solutions in each interval make such a form vanish. This leads to the construction of an aforementioned vanishing equation. These steps also help us to cultivate the effective result on Thue's equations, which address our challenges mentioned earlier.

\section{Integral points on Mordell curves and the proof of Theorem \ref{thm:general-case}}
\label{sec:Mordell-curves-counting}

As discussed previously, we will be employing the counts of integral solutions to Mordell curves in proving Theorem \ref{thm:general-case}.

\subsection{Theorem \ref{thm:general-case} and counting integral points on Mordell curves} In this subsection, we shall prove Theorem \ref{thm:general-case}. Firstly, we state the conjecture on counting integral solutions to Mordell curves. This is what we refer to when we say that we do not have as strong a result as the result for Thue's equation.

\begin{conjecture}[Mordell Curve Conjecture]
\label{conj:weak-effective-Mordell}
For a fixed integer $D\not =0$ and an arbitrarily small $\varepsilon>0$, we have that
\begin{align*}
    \#\{(x,y)\in\Z^2: y^2=x^3+D\} \ll_\eps |D|^{\varepsilon}.
\end{align*}
\end{conjecture}

Using the Helfgott-Vankatesh mechanism \cite{Helfgott-Venkatesh} with the Ellenberg-Venkatesh result \cite{Ellenberg-Venkatesh}, we prove the following proposition, which seems to be the current best record to the best of our knowledge.

\begin{proposition}
\label{prop:upper-bound-Moedell-theorem}
Let $\eps>0$ be given, and let $\varpi_0=\frac{2\beta}{3\log 3}\approx 0.1688$ be as in \eqref{eq:varpi_0-def}. Fixed a nonzero integer $D$, we have that the number of integer solutions to the Mordell curve $y^2=x^3+D$ is at most
$$O_\eps(|D|^{\varpi_0+\eps}).$$
\end{proposition}

\begin{proof}
    We plug the bound of the $3$-torsion of the class group $h_3(\Q(\sqrt{-D}))\ll_\eps D^{\frac{1}{3}+\eps}$ from Proposition 2 of \cite{Ellenberg-Venkatesh} into the mechanism as in the proof of Theorem 4.5 of \cite{Helfgott-Venkatesh}.
\end{proof}

\noindent
\begin{remark}
Proposition \ref{prop:upper-bound-Moedell-theorem} gives a better bound than the significantly improved result on the $2$-torsion of the class group of Bhargava et al. \cite{effective-Mordell}. This is because Theorem 1.2(c) in \cite{effective-Mordell} gives the bound on general elliptic curves rather than Mordell curves. It is interesting if this bound on the $2$-torsion can be exploited to improve upon Conjecture \ref{conj:weak-effective-Mordell}, which might be of its own interest.
\end{remark}

\begin{remark}
It is possible to prove the weak version of Conjecture \ref{conj:weak-effective-Mordell} for special $D$ elementarily. Here, the weak version means that for each $\eps>0$ and a fixed integer $D\not=0$, we have that 
\begin{align}
\label{eq:weak-Mordell-bound}
\#\{(x,y)\in\Z^2: |x|,|y|\leqslant B, y^2=x^3+D\}\ll_{\eps} (B|D|)^{\eps}.
\end{align} 
For any $D$, this weak version can be proved by using Heath-Brown's result as was mentioned in Theorem 2 in \cite{Heath-Brown-counting} in which he proved for non-singular cubic bivariate forms. We shall consider $G(X,Y,Z)=ZY^2-X^3-DZ^3$ in Heath-Brown's theorem. This cubic form covers all solutions of our Mordell curve, and hence, we obtain what we desire. Note that he assumed the Birch-Swinnerton-Dyer conjecture and the Riemann hypothesis for $L$-functions of elliptic curves to obtain the result. Also, this version, and hence invoking two such conjectures, suffices to prove our second assertion in Theorem \ref{thm:general-case}. This is because we may take $B=N$ in the weak version above and use the bound as in the proof below.
\end{remark}

Now, we prove Theorem \ref{thm:general-case}. As in \eqref{eq:split-into-dyadic-of-counting-sqfull}, we consider the majorant $g(x)$ of $f(x)$, say $g(x)=ax^2+b$ with $b\not =0$. Then, we find an upper bound for $\mathcal{M}(E,D)$ by fixing $e\in [E,2E]$ and count the number of solutions $(d,n)$ to $an^2=e^2d^3-b$. By Proposition \ref{prop:upper-bound-Moedell-theorem} with $(x,y)=(e^2ad, e^2a^2n)$ and $D=-e^4ba^3$, we obtain
\begin{align*}
\mathcal{M}(E,D)\ll_{\eps, f} E\cdot |E|^{4\varpi_0+\eps}=E^{1+4\varpi_0+\eps}.
\end{align*}
Thus, we obtain the total bound to be 
$$\ll_{\eps, f} N^\eps\left((\log N)\sup_{\substack{E\leqslant N\\ E^2D^3\asymp N^2}}\min(D,E^{1+4\varpi_0})\right)\ll N^{\varpi+\eps},$$
where $\varpi=\frac{2(1+4\varpi_0)}{5+12\varpi_0}$ for the optimal choice $E=N^{\frac{2}{5+12\varpi_0}}$. The second assertion can be proved similarly by replacing $\varpi_0$ by $0$. The theorem follows as desired.

\noindent
\begin{remark}
Another approach appealing to Heath-Brown's \cite{Heath-Brown-counting} or Salberger's mechanism (as mentioned in \cite{Power-free-values-of-polynomials}) seems unfortunate to tackle this problem with our naive method. We estimate the number of points $(d,n)$ that satisfies
\begin{align*}
q(d)-p(n)-b:=e^2d^3-an^2-b=0,
\end{align*}
where $b\not=0$, $p(x)$ is a polynomial of degree $2$ and $q(x)$ is a polynomial of degree $3$. Then, we proceed as in Reuss's work \cite{consec-square-full} and obtain via theorem 15 of \cite{Heath-Brown-counting} that the points $(d,n)$ satisfying the above lie on at most $N^\eta \sqrt{D}$ auxiliary curves. Thus, using B\'ezout's theorem, the number of points considered is $\ll_\eta N^\eta\sqrt{D}$. Upon invoking Estermann's bound in the first sum, we obtain the upper bound to be $$\ll_\eta N^\eta\left((\log N)\sup_{\substack{E\leqslant N\\ E^2D^3\asymp N^2}}\min(D,E\sqrt{D})\right)\ll N^{\frac{1}{2}+\eta},$$ where the supremum occurs at $E=N^\frac{1}{4}$. This is, however, weaker than what we expected in Theorem \ref{thm:general-case}.
\end{remark} 

\subsection{Conditional result on the $abc$ conjecture}

In this subsection, we show that the number $\varpi$ in Theorems \ref{thm:squarefull-on-quadratic-forms} and \ref{thm:general-case} can be taken arbitrarily close to $0$ under the infamous $abc$ conjecture, stated below.

\begin{conjecture}[$abc$ conjecture]
\label{conj:abc-conjecture}
    For any fixed $\eps>0$ if $a+b=c$ for pairwise coprime positive integers $a,b,c$, then for some constant $K_\eps$
    $$c<K_\eps\cdot\textup{rad}(abc)^{1+\eps}.$$
\end{conjecture}

\noindent

\begin{proposition}
\label{prop:general-polyn-via-abc}
    Let $\eps>0$ and $f(x)$ be an admissible quadratic polynomial. Under the assumption on the $abc$ conjecture, we have that the number of $n\leqslant N$ for which $f(n)$ is square-full is at most $$O_{\eps, f}(N^\eps).$$
\end{proposition}

\begin{proof} 
    Assume that the $abc$ conjecture holds. We may focus on counting the solutions $(n,e,d)$ in $ce^2d^3=n^2+b$ for fixed constants $c,b$ as in Equation \eqref{eq:majorant-g}. Here, we remove the factor $a$ from both $p$ and $q$ in \eqref{eq:majorant-g} and the remaining denominator $4a$ is moved and relabeled as $c$. The sizes of $c$ and $b$ are dependent on $f$. Then, we determine the condition on $d$ that satisfies the equation $ce^2d^3=n^2+b$. For $b>0$, suppose $\gcd(b,n)=\ell$ we transform the equation into
    \begin{align*}
        \frac{ce^2d^3}{\ell\ell_1} = \ell' n'^2+b'',
    \end{align*}
    where $n=\ell n'$ and $b=\ell b'$ with $\gcd(n',b')=1$, and we also let $\ell_1=\gcd(\ell,b')$ with $b'=\ell_1 b'', \ell=\ell_1\ell '$. Thus, there exists an absolute constant $C_\varepsilon$ such that
    \begin{align*}
    \frac{ce^2d^3}{\ell\ell_1}&\leqslant C_\eps\text{rad}\left(\frac{ce^2d^3}{\ell\ell_1}\cdot \ell'n'^2\cdot b''\right)^{1+\eps}
    \\ &\leqslant C_\eps \text{rad}\left(\frac{ce^2d^3}{\ell_1^2}n'b''\right)^{1+\eps}\leqslant C_\eps \text{rad}\left(ce^2d^3n'b''\right)^{1+\eps}\ll_{\eps, f} e^{1+\eps}d^{1+\eps}n'^{1+\eps}.
    \end{align*} 
    Hence, we have $e^{1-\eps}d^{2-\eps}\ll_{\eps, f} n^{1+\eps},$
    and since $e^2d^3\asymp_{f} n^2$, we obtain that $d\ll_{\eps, f} N^\eps$ for arbitrarily small $\eps>0$. Then, the Estermann's bound and the number of possibilities of $d$ imply the claim. For $b<0$, the argument is similar, and we are done.
\end{proof}

\begin{remark}
\label{rm:abc-conj}
Since certain $f(x)$ produce infinitely many square-full numbers, this bound in Proposition \ref{prop:general-polyn-via-abc} is best possible as a power of $N$.
\end{remark}

\section{The determinant method and the proof of Theorem \ref{thm:squarefull-on-quadratic-forms}}
\label{sec:determinant-method-and-special-polyn}

In this section, we prove Theorem \ref{thm:squarefull-on-quadratic-forms}. We proceed as in Heath-Brown \cite{square-free} and Reuss \cite{consec-square-full}. Fixing $\eta>0$, we determine the solutions $(e,d,n)$ to equation $e^2d^3=n^2+\alpha^2$ with $d$ square-free. Recall that we have parameters $E,D,N$ and the conditions $E<e\leqslant 2E, D<d\leqslant 2D$ and $N<n\leqslant 2N$. These provide us $E^2D^3\asymp N^2$ since $\alpha\leqslant N^\varpi$ and $\varpi=\frac{29}{100}$. We have $e=x_1^2+x_2^2$ and $d=y_1^2+y_2^2$ by the sum of two squares argument, and $d$ is square-free. By the unique factorization property of $\Z[i]$, we have that $(x_1+ix_2)^2(y_1+iy_2)^3=n+\alpha i$ and that after taking its imaginary part,
\begin{equation}
\label{eq:imaginary-part-original}
(x_1^2-x_2^2)(3y_1^2y_2-y_2^3)+2x_1x_2(y_1^3-3y_1y_2^2)=\alpha.
\end{equation}
We shall assume that $|x_1|\leqslant |x_2|$ and $x_1,x_2$ have the same signs. For if $|x_1|\geqslant |x_2|$ we may change the role of them and change the sign of $y_2$. If $x_1,x_2$ have different signs, we can change the signs of $y_1$ and $x_2$. Now, we claim that $$\max\left(|3y_1^2y_2-y_2^3|,|y_1^3-3y_1y_2^2|\right)\gg D^{\frac{3}{2}}.$$

The proof will assume $|y_1|\leqslant |y_2|$; however, since they are symmetry the other case can be done in a similar way. For the case $|y_2|\cdot|3y_1^2-y_2^2|\geqslant |y_1|\cdot |y_1^2-3y_2^2|$, we first suppose that $y_2^2<3y_1^2$. Then we have $|y_2|(3y_1^2-y_2^2)\geqslant |y_1||y_1^2-3y_2^2|\geqslant |y_1|(3y_2^2-y_1^2)$.
Since $d$ is square-free we see that $|y_1|\not =|y_2|$, for otherwise $d=2y_1^2$ is not square-free and $d\gg D$ is large for $N\gg_\eta 1$, then $3t^2-1\geqslant t(3-t^2)$ for $1\geqslant t:=|y_1/y_2|$. Thus, we have $0\leqslant t^3+3t^2-3t-1\leqslant 3t^2-3t=3t(t-1)$ and $t=1$ so $|y_1|=|y_2|$, which is not possible. Therefore, $y_2^2\geqslant 3y_1^2$ and with the same setup we obtain that $0\leqslant t^3-3t^2-3t+1=(t-1)^3-6t+2\leqslant 2-6t$ and $|y_1|\leqslant \frac{1}{3}|y_2|$. The maximum in this case is $\gg |y_2|(y_2^2-3y_1^2)\gg D^{3/2}$. Similarly, for the case $|y_2||3y_1^2-y_2^2|<|y_1||y_1^2-3y_2^2|$, with the same notation we have either
$$t(3-t^2)>3t^2-1\text{ or } t(3-t^2)>1-3t^2,$$
and the second inequality is impossible, as it is equivalent to $(t+1)(t^2-4t+1)<0$. The first corresponds to $(t-1)(t^2+4t+1)<0$, which is always true, and corresponds to $t^2<3\Longrightarrow y_2^2<3y_1^2$. Again, we obtain the maximum of the two terms $\gg D^{3/2}$.
\\

Upon taking $q_1(y_1,y_2)$ to be the largest of the two, which has the property that its absolute value is $\gg D^{3/2}$, and taking $z_1,z_2$ appropriately to be polynomials in $x_1,x_2$, we obtain
\begin{align*}
    q_1(y_1,y_2)z_1+q_2(y_1,y_2)z_2=\alpha.
\end{align*}

We may assume that $D\geqslant N^\varpi$. Estermann \cite{square-free-Estermann} provides the bound on
$$\#\{(n,e)\in\Z^2:|n|,|e|\leqslant N, n^2=d^3e^2-\alpha^2\}\ll \log (N+|d\alpha|),$$
which is his Equation (16). For $D<N^\varpi$, it contributes
\begin{align*}
    \sum_{\substack{N<n\leqslant 2N\\ e^2d^3=n^2+\alpha^2\\ D<d\leqslant 2D}}\mu^2(d)\leqslant\sum_{d\leqslant 2N^\varpi} \sum_{\substack{n^2=d^3e^2-\alpha^2\\|n|,|e|\leqslant 2N}} 1\ll N^\varpi\log N,
\end{align*}
which is acceptable. By the triangle inequality, we have that $$D^{\frac{3}{2}}|z_1|\ll |q_1(y_1,y_2)||z_1|\leqslant\alpha +|q_2(y_1,y_2)||z_2|\ll D^{\frac{3}{2}}|z_2|,$$
since $\alpha\leqslant N^\varpi\leqslant D$ and here we assume that $|z_1z_2|\not= 0$. For if it is the case, then $x_1=x_2$, and \eqref{eq:imaginary-part-original} becomes $2x_1^2y_1(y_1^2-3y_2^2)=\alpha$. This implies that the number of $(x_1,x_2,y_1,y_2)$ satisfying such an equation is $\ll\tau(\alpha)^{10}\ll N^\eps$, which contributes in an acceptable bound of counting $(e,d,n)$. Thus, $|z_1|\ll |z_2|$ and $|z_2|\gg\max(|z_1|,|z_2|)\gg E$. Hence,
\begin{equation}
\label{eq:imaginary-part-equation}
t=-\dfrac{q_2(s,1)}{q_1(s,1)}+O(E^{-1}D^{-\frac{3}{2}})=-\dfrac{q_2(s,1)}{q_1(s,1)}+O(N^{-1})=:\phi(s)+O(N^{-1}),
\end{equation}
where $s=y_1/y_2$ and $t=z_1/z_2$ if $|y_1|\leqslant |y_2|$ and $s=y_2/y_1$ in the other case, which can be handled in a similar way, so we may assume $s,t\ll 1$. Therefore, the problem reduces to counting points $(s,t)$ that lie close to the curve $t=\phi(s)$. We write $w=x_1/x_2$ which is positive since $x_1,x_2$ have the same signs and $w\ll 1$.

\subsection{The determinant method}

In this subsection, we construct another equation in a bi-homogeneous form of the same variables as in \eqref{eq:imaginary-part-equation} for each interval that we will be splitting. The satisfactory bound is plausible since we have the existing bound for counting points on $G(x_1,x_2;y_1,y_2)=0$, for an absolutely irreducible bi-homogeneous polynomial $G(x_1,x_2;y_1,y_2)\in\Z[x_1,x_2,y_1,y_2]$.\\

Our plan is as follows. We proceed to choose a parameter $M\leqslant N$, which will be specified later in terms of $E,D$. Then, we divide the intervals of possible values of $s\ll 1$ into $O(M)$ intervals of the form $I=(s_0, s_0+M^{-1}]$. By Taylor's theorem, if we write $s\in I$ by $s=s_0+u$ with $u=O(M^{-1})$, then $\phi(s)=\phi(s_0)+P(u)+O(N^{-1})$ where $P(u)$ is a polynomial in $u$ without a constant term and all coefficients of size $O(1)$. Thus, $|P(u)|\ll M^{-1}$. Write $\tau(w):=t=\phi(s)+O(N^{-1})$ where $0<w=x_1/x_2\ll 1$ from our conditions. If $\tau(w)=(w^2-1)/2w$ then we obtain
$$\dfrac{w^2-1}{2w}=\phi(s_0)+P(u)+O(N^{-1}).$$
Since $w>0$ and $P(u)$ is small, we have that
$$w=\phi(s_0)+P(u)+\sqrt{(\phi(s_0)+P(u))^2+1+O(N^{-1})}.$$
By Taylor's theorem, $\sqrt{s+c}=\sqrt{c}+(2\sqrt{c})^{-1}s+O(s^2)$ for $s\in\R$ around zero and $c\geqslant 1$ is a constant. Hence,
\begin{align*}
w &= \phi(s_0)+P(u)+\sqrt{1+(\phi(s_0)+P(u))^2}+\dfrac{O(N^{-1})}{2\sqrt{1+(\phi(s_0)+P(u))^2}}
\\ &= c_{s_0}+Q(u)+O(N^{-1}),
\end{align*}
where we have used Taylor's theorem for the square root above, and the equation holds for some polynomial $Q$ in $u$. Similarly, if $\tau(w)=2w/(w^2-1)$, we obtain
\begin{align*}
w=\dfrac{1+\sqrt{1+(\phi(s_0)+P(u)+O(N^{-1}))^2}}{\phi(s_0)+P(u)}\ll 1,
\end{align*}
and as $P(u)$ is small, $|\phi(s_0)|\gg 1$. By Taylor's theorem on $\sqrt{1+(c+s)^2}$ for $s\in\R$ around zero and $c\geqslant 1$ is a constant, we have
\begin{align*}
w &= \dfrac{1+\sqrt{1+(\phi(s_0)+P(u))^2}}{\phi(s_0)+P(u)}+\dfrac{O(N^{-1})}{\sqrt{1+(\phi(s_0)+P(u))^2}}
\\ &=c_{s_0}+Q(u)+O(N^{-1}),
\end{align*}
by the fact that $|\phi(s_0)|\gg 1$, and we have used Taylor's theorem for the first term on $1+\sqrt{1+(\phi(s_0)+s)^2}$ for $s=P(u)$ is around zero. Since $P(u)\ll M^{-1}$ and $M$ will be chosen so that $\log M\asymp \log N$, we need finite terms accumulated in $Q(u)$ with the residual term of order $O(N^{-1})$. Therefore, in either case, $w$ can be represented by a polynomial $Q(u)$ such that $Q(u)\ll M^{-1}$ with error $O(N^{-1})$.\\

Suppose that there are $J$ solutions to \eqref{eq:imaginary-part-equation} in the interval $I$. Upon introducing new parameters $K,L$ that will be dependent merely on $\eta$, we construct a $J\times H$ matrix $$\mathfrak M=
\begin{pmatrix}
1 & s_1 & w_1 & \cdots & s_1^k w_1^\ell & \cdots & s_1^K w_1^L\\
1 & s_2 & w_2 & \cdots & s_2^k w_2^\ell & \cdots & s_2^K w_2^L\\
\vdots & & & \vdots & & & \vdots \\
1 & s_J & w_J & \cdots & s_J^k w_J^\ell & \cdots & s_J^K w_J^L\\
\end{pmatrix},$$
where $H=(K+1)(L+1)$. We obtain a nonzero vector $\mathbf{c}$ such that
$$\mathfrak M\mathbf{c}=\mathbf{0},$$
if we have $\text{rank}(\mathfrak M)<H$. With this nonzero vector $\mathbf{c}$, we can consider the polynomial $C_I(s,w)=\sum_{h=1}^H c_hm_h(s,w)$, where we label $m_h(s,w)$ to be each monomial $s^kw^\ell$ appearing in each row of $\mathfrak{M}$ for $k\leqslant K$ and $\ell\leqslant L$ so that $C_I(s_j, w_j)=0$ for all $j\leqslant J$. Also, by reducing matrix $\mathfrak{M}$ into row-reduced echelon form, we have that $\mathbf{c}\in\Q^H$ has rational entries and we may clean up the denominator of $C_I$, which is $\ll$ some fixed power of $N$, so that its coefficients are integral of size $\ll N^{\kappa(K,L)}$. The fact that $\text{rank}(\mathfrak{M})$ is strictly less than $H$ follows from $$\left(\prod_{j\leqslant H} y^K_{2,j}x^L_{2,j}\right)\Delta\in\Z,$$
where $\Delta$ is the subdeterminant of the first $H\times H$ matrix, and we assume that $J\geqslant H$. Thus, we will choose $M$ appropriately so that $\Delta=0$ for which it suffices to prove that $\Delta\ll_{K,L} D^{-\frac{KH}{2}}E^{-\frac{LH}{2}}$ for the appropriate implied constant, and we obtain the desired inequality. 
\\

Following the proof of Heath-Brown \cite{square-free}, for $V\asymp N$ in the Taylor's expansion $w=c_{s_0}+P(u)+v$ with $v\ll V^{-1}$ and note that $u\ll M^{-1}$, we now order the values of $M^{-k}V^{-\ell}$, $k\leqslant K$, $\ell\leqslant L$, decreasingly as $1=M_0,M_1,\dots$. Lemma 3 of Heath-Brown \cite{Heath-Brown-lemma} gives
$$\Delta\ll_H \prod_{h\leqslant H} M_h.$$
Let us denote $M_H=W^{-1}$. Then $M^{-k}V^{-\ell}\geqslant M_H$ if and only if 
$$k\log M+\ell\log V\leqslant \log W.$$
The number of pairs $k,\ell$ satisfying the above is
$$\dfrac{\log^2 W}{2\log M\log V}+O\left(\dfrac{\log W}{\log N}\right)+O(1),$$
which equals $H$. Since $\log M,\log V\ll \log N$, we have that
$$\log W=\sqrt{2H\log M\log V}+O(\log N).$$
Moreover, we see that
\begin{align*}
\log\prod_{h\leqslant H} M_h 
 &=\sum_{\substack{k.\ell \\ k\log M+\ell \log V\leqslant \log W}} -(k\log M+\ell\log V)
\\ &= -\dfrac{\log^3 W}{3\log M\log V}+O\left(\dfrac{\log^2 W}{\log N}\right)
\\ &=-H^{\frac{3}{2}}\cdot\dfrac{2\sqrt{2}}{3}\sqrt{\log M\log V} +O(H\log N),
\end{align*}
by replacing the term of $\log W$ above. Whence, we acquire 
$$\log |\Delta|\leqslant O_H(1)-H^{\frac{3}{2}}\cdot\dfrac{2\sqrt{2}}{3}\sqrt{\log M\log V}+O(H\log N).$$
Thus, $\Delta=0$ if we have that \begin{align*}
\dfrac{K}{2}\log D+\dfrac{L}{2}\log E\leqslant (KL)^{\frac{1}{2}}\dfrac{2\sqrt{2}}{3}\sqrt{\log M\log V}+O_{K,L}(1)+O(\log N),
\end{align*}
or, for taking $K=[L\log E/\log D]$ that
\begin{align*}
L\log E\leqslant L\cdot\dfrac{2\sqrt{2}}{3}\sqrt{\log M\log V}\sqrt{\dfrac{\log E}{\log D}}+O_{L}(1)+O(\log N).
\end{align*}

For small fixed $\delta>0$ if
$$\dfrac{2\sqrt{2}}{3}\sqrt{\log M\log V}\sqrt{\dfrac{\log E}{\log D}}\geqslant (1+\delta)\log E,$$
and $L,N\gg_\delta 1$, we will obtain the desired inequality. Since $V\gg N$, we summarize and acquire the following lemma.
\begin{lemma}
\label{lem:M-and-helping-equation}
Let $\eta>0$ be given and $M\in\left[D, N\right]$ satisfies
$$\log M\geqslant \dfrac{9}{8}(1+\eta)\dfrac{\log E\log D}{\log N}.$$
Then, for any interval $I=[s_0,s_0+M^{-1}]$, there exists $0\not =C_I(s,t)$ of integral coefficients such that all solutions satisfy
$$C_I\left(y_1/y_2,x_1/x_2\right)=0,$$
with $y_1/y_2\in I$. Moreover, $\deg(C_I)=O_\eta(1)$ and its coefficients size $O_\eta(N^{\kappa})$ for some $\kappa=\kappa(\eta)$.
\end{lemma}

With this polynomial $C_I$ in the previous lemma, as in \cite{square-free} we may assume without loss of generality that $C_I$ is absolutely irreducible. Hence, we can apply the existing bound on counting solutions noted at the beginning of this subsection.

\subsection{Counting the solutions and finishing the proof}
In the previous subsection, we construct an auxiliary equation involving a bi-homogeneous polynomial $C_I$ in each fixed interval $I$. It can be proven that we may restrict to $C_I$ that is absolutely irreducible, say $F(y_1,y_2;z_1,z_2)=0$ to be again bi-homogeneous. In this subsection, we will count the solutions to \eqref{eq:imaginary-part-original} by changing of bases and employing such a constructed vanishing equation. We now have in the interval $I=(s_0,s_0+M^{-1}]$ that
$$s_0<s=\dfrac{y_1}{y_2}\leqslant s_0+M^{-1}.$$
It follows that $|y_1-s_0y_2|\leqslant D^{1/2}M^{-1}$ because $|y_2|\leqslant D^{1/2}$. We proceed by letting 
$$\Lambda = \left\{(D^{-1/2}M(y_1-s_0y_2), D^{-1/2}y_2):(y_1,y_2)\in\Z^2\right\}.$$
Then $\Lambda$ is a lattice of determinant $D^{-1}M$. We are interested in points $(\alpha_1,\alpha_2)\in\Lambda$ falling in the square
$$S=\{(\alpha_1,\alpha_2):\max(|\alpha_1|,|\alpha_2|)\leqslant 1\}.$$

Let $\mathbf{g}^{(1)}$ be the shortest non-zero vector in the lattice and $\mathbf{g}^{(2)}$ be the shortest vector not parallel to $\mathbf{g}^{(1)}$. The two vectors form a basis for $\Lambda$ and, moreover, $\lambda_1\mathbf{g}^{(1)}+\lambda_2\mathbf{g}^{(2)}\in S$ only when $|\lambda_1|\ll |\mathbf{g}^{(1)}|^{-1}=: L_1$ and $|\lambda_2|\ll|\vecg^{(2)}|^{-1}=:L_2$. Thus, we have $L_1\gg L_2$ and $(L_1L_2)^{-1}=|\vecg^{(1)}||\vecg^{(2)}|\ll\det(\Lambda)=D^{-1}M\Longrightarrow L_1L_2\gg DM^{-1}$. Hence, we may replace $(y_1,y_2)$ by $(\lambda_1,\lambda_2)$. We may argue in exactly the same way to change $x_1,x_2$ by $\tau_1,\tau_2$. Now, we have the equations
\begin{align}
\label{eq:equal-1}
G_0(\lambda_1,\lambda_2;\tau_1,\tau_2) &= \alpha \\
\label{eq:equal-0}
G_1(\lambda_1,\lambda_2;\tau_1,\tau_2) &= 0
\end{align}
which correspond, respectively, to \eqref{eq:imaginary-part-original} and the vanishing equation from Lemma \ref{lem:M-and-helping-equation}. $G_0,G_1$ are bi-homogeneous polynomials of degree $(3;2)$ and $(a;b)$, respectively. When $a,b\geqslant 2$ we can get a satisfactory bound from the following lemma due to Lemma 2 of Heath-Brown \cite{square-free}.
\begin{lemma}
\label{lem:number-of-sols-equal-0}
Let $G(x_1,x_2;y_1,y_2)$ be an absolutely irreducible bi-homogeneous polynomial with integer coefficients of degree $(a; b)$ with $a,b\geqslant 1$. Let $\varepsilon>0$ be given. Then for $X\geqslant 1$ there exist $O_{a,b,\varepsilon}(X^{\frac{2}{b}+\varepsilon}\lVert G\rVert^\varepsilon)$ points $(a_1,a_2,b_1,b_2)\in\Z^4$ satisfying the conditions
\begin{align*}
\gcd(a_1,a_2) &= \gcd(b_1,b_2)=1,
\\ G(a_1,a_2;b_1,b_2) &= 0 \text{ and } \max\{a_1,a_2\}\leqslant X.
\end{align*}
\end{lemma}

It is possible to assume that $\gcd(\lambda_1,\lambda_2)=\gcd(\tau_1,\tau_2)=1$ at the expense of the factor of $\alpha^\eps$. Suppose that $\gcd(\lambda_1,\lambda_2)=\lambda$ then $\lambda^3|\alpha$ as $G_0$ is homogeneous in $\lambda_1,\lambda_2$. If this happens with $\lambda>1$ we split this case and count the solutions from it. There are $O(\tau(\alpha))=O(\alpha^{\eps/2})$ possibilities of such equations. Similarly, there exist $O(\alpha^{\eps/2})$ possible equations when we determine a pair of $\tau_1,\tau_2$.\\

We present here an elementary result on counting integral solutions with some conditions on certain bivariate forms.
\begin{lemma}
\label{lem:irr-y^3}
Let $N$ be a large positive number and $\alpha$ be as in Theorem \ref{thm:squarefull-on-quadratic-forms}. Let $c,d$ be a pair of integers with $|c|,|d|\leqslant N^{1000}$. Then, the number of solutions $(y_1,y_2)\in\Z^2$ such that $|y_1|,|y_2|\leqslant N$ to the equation 
\begin{align*}
    P_{c,d}(y_1,y_2):=c(3y_1^2y_2-y_2^3)+d(y^3_1-3y_1y_2^2)=\alpha
\end{align*}
is $O_\eps(N^\eps)$. 
\end{lemma}
\begin{proof}
Without loss of generality, we assume that $\gcd(c,d)=1$ since $\alpha\leqslant N^{\frac{29}{100}}$ and there are at most $\tau(\alpha)\ll_\eps N^\eps$ possible cases to determine. We first prove that $P(y_1,y_2):=P_{c,d}(y_1,y_2)$ is not a product of linear factors. Suppose that there exists a linear factor $qy_1-py_2$ with $p,q\in\Z$ and $\gcd(p,q)=1$. Then, we have that
    \begin{align*}
        dp^3+3cp^2q-3dpq^2-cq^3=0
    \end{align*}
    which implies that $p|c$ and $q|d$ so let $c=pc',d=qd'$ with $\gcd(c',d')=1$. By plugging this in the equation, we obtain
    \begin{align*}
        \dfrac{c'+3d'}{3c'+d'}=\frac{p^2}{q^2},
    \end{align*}
    for which by parity checking, we have that $\gcd(c'+3d',3c'+d')=1,2$ or $4$. If we restrict ourselves to the cases $\gcd(c'+3d',3c'+d')=1,2$, we have either $c'+3d'=p'^2$ and $3c'+d'=q'^2$ or $c'+3d'=2p'^2$ and $3c'+d'=2q'^2$ are our possible systems of equations in integers. Then, we have in each case that $p'^2-3q'^2\equiv 0\pmod{4}$, which is impossible for $\gcd(p',q')=1$. In this case $P(y_1,y_2)$ is irreducible over $\Q$ and we obtain from the main theorem of Bombieri-Schmidt \cite{bombieri-schmidt} involving Thue's equation that the bound is $O(3^{1+\sup_{d|\alpha}\omega(d)})$, where $\omega(d)$ is the number of distinct prime factors of $d$. We have that $3^{\sup_{d|\alpha} \omega(d)}\leqslant 3^{\omega(\alpha)}\leqslant \tau(\alpha)^{\frac{\log 3}{\log 2}}\ll_\eps N^\eps$. The total bound is therefore $O_\eps(N^{\eps})$.
    
    When $\gcd(c'+3d',3c'+d')=4$, we have that
    \begin{align*}
        c = \frac{p(3q^2-p^2)}{2}, d = \frac{q(3p^2-q^2)}{2}.
    \end{align*}
    Then, we use the above to compute the discriminant of the quadratic polynomial
    \begin{align*}
        \frac{P(y_1,y_2)}{qy_1-py_2}=:ay_1^2+by_1y_2+ey_2^2\in\Q[y_1,y_2],
    \end{align*}
    and find that $b^2-4ae=3(p^2-q^2)^2$. This means that it is not possible to factor further into $2$ linear factors over $\Q$. On the flip side, if there is no such linear factor, $P(y_1,y_2)$ is irreducible over $\Q$. Thus, we showed that $P(y_1,y_2)$ is not a product of linear factors.\\

    Now, we are left to consider when there is a quadratic factor of $P(y_1,y_2)$. In this case, we bound the number of solutions of $ay_1^2+by_1y_2+ey_2^2=\beta$ for a certain divisor $\beta$ of $\alpha$. Since there are at most $\tau(\alpha)\ll_\eps N^{\eps}$, we may consider only a fixed case $\beta|\alpha$. Such an equation can be rearranged into
    \begin{align}
    \label{eq:counting-P_{c,d}}
        (2ay_1+by_2)^2-(4ae-b^2)y_2^2=4a\beta.
    \end{align}
    The number of pairs $(y_1,y_2)$ satisfying \eqref{eq:counting-P_{c,d}} is at most the number of pairs $(z_1,z_2)\in \Z^2$ satisfying $z_1^2-\gamma z_2^2:=z_1^2-(4ae-b^2)y_2^2=4a\beta$ and $|z_1|,|z_2|\ll N$. We obtain the bound $O_\eps(N^\eps)$ in this case since $|\gamma|\ll N^{2000}$ and $4a\beta\ll N^{1001}$. Hence, the lemma follows.
\end{proof}

Now, let us count the number of solutions to Equations \eqref{eq:equal-1} and \eqref{eq:equal-0}. When $b\geqslant 2$, we set $X=L_1$ in Lemma \ref{lem:number-of-sols-equal-0}. In Equation \eqref{eq:equal-0}, Lemma \ref{lem:M-and-helping-equation} implies $\lVert G_1\rVert\ll N^{\kappa}$ where $\kappa=\kappa(\eta)$. Then, Lemma \ref{lem:number-of-sols-equal-0} tells us that the number of solutions to \eqref{eq:equal-0} is $O_\eta(L_1^{1+\eta}N^\eta)$. Each of these solutions produces at most one solution in \eqref{eq:equal-1}. Therefore, the total number of solutions satisfying both equations is $O_\eta(N^\eta L_1^{1+\eta})$. In a similar way, we obtain for $a\geqslant 2$ that the contribution is $O_\eta(N^\eta T_1^{1+\eta})$.\\

For $b=0$, Equation \eqref{eq:equal-0} becomes $\tilde{G}_1(\lambda_1,\lambda_2)=0$ for some homogeneous polynomial $\tilde G_1$. This corresponds to $O(1)$ pairs of $(\lambda_1, \lambda_2)$, since $\gcd(\lambda_1,\lambda_2)=1$. Each of such $d$ produces $O(\log N)$ pairs of $(n,e)$ satisfying $\frac{N}{2}\leqslant n<N$ and $\frac{E}{2}\leqslant e<E$ in the equation $e^2d^3=n^2+\alpha^2$ from the Estermann bound. The total contribution is $O(\log N)$. In the case $a=0$, we have that Equation \eqref{eq:equal-0} becomes $\tilde{G}_1(\tau_1,\tau_2)=0$ for a certain homogeneous polynomial $\tilde{G}_1$. The number of possible $(\tau_1,\tau_2)$ is $O(1)$ in this case. We then substitute each $(\tau_1,\tau_2)$ into \eqref{eq:equal-1} and consider counting $(y_1,y_2)$ in \eqref{eq:imaginary-part-original} for fixed $(x_1,x_2)$. By Lemma \ref{lem:irr-y^3}, the contribution is hence totally $O_\eta(N^\eta)$.\\

Then, upon considering the cases $b=1$ and $a=1$ as in the work of Heath-Brown (the discussion before Lemma 3 of \cite{square-free}), we obtain the following lemma.
\begin{lemma}
\label{lem:contribution-to-M(E,F)-in-each-interval}
For $\eta>0$, in an interval $I$, $\mathcal M(E,D)\ll_\eta N^\eta \min(L_1^{1+\eta},T_1^{1+\eta})$.
\end{lemma}

Now, we finish the proof of Theorem \ref{thm:squarefull-on-quadratic-forms}. After fixing a parameter $M$ as in Lemma \ref{lem:M-and-helping-equation} we now count possible intervals by fixing $L\ll N$ and consider each $L< L_1\leqslant 2L$. Then we sum the total contribution dyadically on $L$ by applying Lemma \ref{lem:contribution-to-M(E,F)-in-each-interval}. If $(y_1,y_2)$ corresponds to $\vecg^{(1)}$, then $L_1(y_1-s_0y_2)\ll M^{-1}\sqrt{D}$ and $L_1y_2\ll \sqrt{D}$. Suppose that for some fixed $\varepsilon>0$, $L\gg D^{\frac{1}{2}+\varepsilon}$ we have that $y_2=0$ since $L|y_2|\ll L_1|y_2|\ll\sqrt{D}$ and $y_2\in\Z$. Also, we have that $y_1=0$ as well in the first inequality, which is not possible. Whence, $L_1\ll_\varepsilon D^{\frac{1}{2}+\varepsilon}$ for any $\varepsilon>0$. Then, we write $s_0=y_3M^{-1}$ in $I=(s_0,s_0+M^{-1}]$ so that $y_3\ll M$ since $s_0\ll 1$. Therefore, the number of intervals for which $L<L_1\leqslant 2L$ is at most the number of triples $(y_1,y_2,y_3)\in\Z^3$ with $\gcd(y_1,y_2)=1$ for which
$$y_2y_3=My_1+O(L^{-1}\sqrt{D}), y_2\ll L^{-1}\sqrt{D}, \text{ and } y_3\ll M.$$

Recall that $L_1\gg L_2$ and $L_1L_2\gg DM^{-1}$, so $L\gg D^{\frac{1}{2}}M^{-\frac{1}{2}}$. For each value of $y_1$, there are $O_\eta(N^{\eta}L^{-1}\sqrt{D})$ pairs of $(y_2,y_3)$. Hence, for which $L_1$ is of order $L$ there are $O_\eta(N^\eta L^{-2}D)$ intervals since $y_1\ll L^{-1}\sqrt{D}$. We find from Lemma \ref{lem:contribution-to-M(E,F)-in-each-interval} that $\mathcal{M}(E,D)\ll_\eta N^{2\eta}L^{-2}D(L^{1+\eta})\ll N^{3\eta}\sqrt{DM}$ since $L\gg D^{\frac{1}{2}}M^{-\frac{1}{2}}$. Analogously, we obtain $\mathcal{M}(E,D)\ll N^{3\eta}\sqrt{ME}$. Therefore, we conclude that $\mathcal{M}(E,D)\ll N^{3\eta}\sqrt{M}\min(\sqrt{E},\sqrt{D}).$\\

Next, we let $M$ satisfy the condition in Lemma \ref{lem:M-and-helping-equation}, i.e., dependent on $E,D$. By the consequence of Lemma \ref{lem:contribution-to-M(E,F)-in-each-interval} on $\mathcal{M}(E,D)$ above, we have that the total contribution in each dyadic range is at most
\begin{align*}
\dfrac{1}{N^\eta}&\left(\sum_{\ell\ll \log N}\sum_{\substack{\frac{N}{2^{\ell}}<E\leqslant \frac{N}{2^{\ell-1}}\\E^2D^3\asymp N^2}} \sqrt{M}\min(\sqrt{E},\sqrt{D})\right)
\\ &\ll_\eta \dfrac{1}{N^\eta}\left((\log N)\sup_{\substack{\psi\in[0,1]}} N^{\frac{3}{8}\psi(1-\psi)+\min(\frac{\psi}{2},\frac{1-\psi}{3})}\right)\ll_\eta N^{\varpi},
\end{align*}
where $\varpi=\frac{29}{100}$ for the optimal choice $\psi=\frac{2}{5}$. This completes the proof.
\noindent
\begin{remark}
\label{rm:rand-polynom}
The $abc$ conjecture implies $\sqfull(N)\ll_{\eps,f} N^\eps$; however, we hope for the asymptotic relation, which is much smaller than that. It is possible to see the behavior of average quadratic polynomials, so we determine square-full values of random polynomials. Here, we consider ordering polynomials via naive height similar to the work of Browning and Shparlinski \cite{square-free-random}. For positive integers $H$ and $k\geqslant 2$, we write $\textbf{a}=(a_0,a_1,\dots,a_k)\in\Z^{k+1}$ and denote
$$\mathcal{F}_k(H):=\{a_0+a_1X+\cdots+a_kX^k\in\Z[X]: \textbf{a}\in\mathcal{B}_k(H)\},$$
where $\mathcal{B}_k(H):=\{\mathbf{a}\in\Z^{k+1}:|a_i|\leqslant H,\text{ for all } i=0,1,\dots,k\}$. Here, we focus only on $k=2$. Then, we have by a naive bound fixing $a_1,a_2,n$ and determining the square in the interval
\begin{align*}
    e^2=\frac{a_0+a_1n+a_2n^2}{d^3}\in\left[\frac{a_1n+a_2n^2}{d^3}, \frac{a_1n+a_2n^2}{d^3}+\frac{H}{d^3}\right].
\end{align*}
The number of squares in such an interval is $O(\sqrt{Hd^{-3}}+1)$ and each $e$ produces at most $1$ solution of $a_0$. Whence, we obtain that
\begin{align*}
\sum_{f\in\mathcal{F}_2(H)}\sqfull(N) &\ll \sum_{d\ll N^{\frac{2}{3}}}\#\{(e,d,n,\textbf{a}):n\leqslant N, \dfrac{a_0+a_1n+a_2n^2}{d^3}=e^2\}
\\ &\ll \sum_{d\ll N^{\frac{2}{3}}} \left(\sqrt{\dfrac{H}{d^3}}+1\right)H^2N \ll H^{\frac{5}{2}}N+H^2N^{\frac{5}{3}}.
\end{align*}
Hence, we obtain for small $\eps>0$ that whenever $N^{2+\eps}\leqslant H\leqslant N^A$ for some constant $A$,
$$\dfrac{1}{\#\mathcal{F}_2(H)}\sum_{f\in\mathcal{F}_2(H)}\sqfull(N)\ll_\eps \dfrac{1}{N^\frac{\eps}{2}}.$$
This means that almost all quadratic polynomials have pretty small bounds for their counting formulae. 
\end{remark}

\subsection*{Acknowledgement}
The authors are grateful to Professor Yuk-Kum Lau and Professor Ben Kane for encouraging discussions and the organizers for their hospitality during the HKU Number Theory Days 2024. We would also like to thank Professor Bryce Kerr for interesting discussions during the presentation of this work, as well as the organizers of the Number Theory Down Under 12 conference for their warm hospitality.

\printbibliography[
heading=bibintoc,
title={References}
]

@misc{square-free,
    key="Heath-Brown",
    note={D. R. Heath-Brown, Square-free values of $n^2+1$, \textit{Acta Arith.} \textbf{155} (2012) 1--13.}
}

@misc{Heath-Brown-counting,
    key="Heath-Brown",
    note = {D. R. Heath-Brown, \textit{Counting Rational Points on Algebraic Varieties}, Springer-Verlag, 2006}
}

@misc{Bateman-Grosswald,
    key="Bateman",
    note={P. T. Bateman and E. Grosswald, On a theorem of Erd\H{o}s and Szekeres, \textit{Illinois J. Math.} \textbf{2} (1958) 88--98.}
}

@misc{square-full-inAP,
    key="Chan",
    note={T. H. Chan, Squarefull numbers in arithmetic progression II, \textit{J. Number Theory} \textbf{152} (2015) 90--104.}
}

@misc{square-free-Estermann,
    key="Estermann",
    note={T. Estermann, Einige Sätze über quadratfreie Zahlen (German), \textit{Math. Ann.} \textbf{105} (1931) 653--662.}
}

@misc{square-full-Erdos-Szekeres,
    key="Erd\H{o}s",
    note={P. Erd\H{o}s and G. Szekeres, \"{U}ber die Anzahl der Abelschen Gruppen gegebener Ordnung und \"{u}ber ein verwandtes zahlentheoretisches Problem (German), \textit{Acta Sci. Math. (Szeged)} \textbf{7} (1934--1935) 95--102.}
}

@misc{effective-Mordell,
    key="Bhargava",
    note={M. Bhargava and A. Shankar, et al., Bounds on $2$-torsion in class groups of number fields and integral points on elliptic curves, \textit{J. Am. Math. Soc.} \textbf{33} (2017) 1087--1099.}
}

@misc{consec-square-full,
    key="Reuss",
    note={T. Reuss, Pairs of $k$-free numbers, consecutive square-full numbers, \textit{Arxiv.} (2014).}
}

@misc{Heath-Brown-lemma,
    key="Heath-Brown",
    note={D. R. Heath-Brown, Sums and differences of three $k$-th powers, \textit{J. Number Theory} \textbf{129} (2009) 1579--1594.}
}

@misc{square-free-random,
    key="Browning",
    note={T. Browning and I. Shparlinski, Square-free values of random polynomials, \textit{J. Number Theory} \textbf{261} (2024) 220--240.}
}

@misc{Power-free-values-of-polynomials,
    key = "Browning",
    note ={T. Browning, Power-free values of polynomials, \textit{Arch. Math.} \textbf{96} (2011) 139--150.} 
}

@misc{Granville,
    key ="Granville" ,
    note = {
        A. Granville, ABC allows us to count squarefrees, \textit{Int. Math. Res.} \textbf{19} (1998)  991--1009
    }
}

@misc{bombieri-schmidt,
    key ="Bombieri" ,
    note = {
        E. Bombieri and W. M. Schmidt, On Thue's equation, \textit{Invent. Math.} \textbf{88} (1987)  69--81
    }
}

@misc{Helfgott-Venkatesh,
    key ="Helfgott" ,
    note = {
        H. A. Helfgott and A. Venkatesh, Integral points on elliptic curves and 3-torsion in class groups, \textit{J. Amer. Math. Soc.} \textbf{19} (2006)  527--550
    }
}

@misc{Ellenberg-Venkatesh,
    key ="Ellenberg" ,
    note = {
        J. Ellenberg and A. Venkatesh, Reflection principles and bounds for class group torsion, \textit{Int. Math. Res. Not.} (2007)
    }
}

\end{document}